\documentclass[letterpaper]{amsart}
\usepackage{amssymb}
\usepackage{amsmath}

\newtheorem{theorem}{Theorem}[section]
\newtheorem{lemma}[theorem]{Lemma}
\newtheorem{proposition}[theorem]{Proposition}
\newtheorem{corollary}[theorem]{Corollary}
\newtheorem*{theorem*}{Theorem}

\theoremstyle{definition}

\newtheorem{assumption}{Assumption}

\newcommand{\euO}{\mathfrak O}
\newcommand{\euP}{\mathfrak P}

\newcommand{\bF}{\mathbb F}

\begin{document}
\title[Galois scaffolding]{One-dimensional elementary abelian extensions have Galois scaffolding}

\author{G. Griffith Elder} 
\email{elder@vt.edu} 
\address{Department of
Mathematics \\ University of Nebraska at Omaha\\ Omaha, NE 68182-0243
U.S.A.}  \curraddr{Department of Mathematics \\ Virginia Tech\\
Blacksburg, VA 24061-0123 U.S.A.}
\date{May 2, 2007}
\subjclass{11S15} \keywords{Ramification, Galois module structure}
\thanks{The author was partially supported by National Science Foundation Grant No. 201080.}

\bibliographystyle{amsalpha}

\begin{abstract}
We define a variant of normal basis, called a {\em Galois
scaffolding}, that allows for an easy determination of valuation, and
has implications for Galois module structure.
We identify fully ramified, elementary abelian extensions of local
function fields of characteristic $p$, called {\em one-dimensional},
that, in a particular sense, are as simple as cyclic degree $p$
extensions, and prove the statement in the title above.
\end{abstract}

\maketitle

\section{Introduction}

The Normal Basis Theorem states that in a finite, Galois extension
$L/K$ with $G=\mbox{Gal}(L/K)$,
there are elements $\rho\in L$ whose conjugates
$\{\sigma\rho:\sigma\in G\}$ provide a field basis for
$L$ over $K$.  In the setting of local field extensions, the most
important property of an element is its valuation and so we asked in
\cite{elder:blms} about the valuation of these elements: Are there are
valuations (integer certificates) that guarantee that any element
bearing the specified valuation be a normal basis generator? ({\em
i.e.}  $v\in \mathbb{Z}$ so that $\rho\in L$ and $v_L(\rho)=v$ implies
$\{\sigma\rho:\sigma\in G\}$ is a basis for $L$ over
$K$.)

In this paper, we ask for more.  Let $L/K$ be a fully ramified
$p$-extension of local fields with finite residue field of
characteristic $p$, and let $v_L$ denote the normalized, additive
valuation.  We ask, in addition to the above property, that there be
an explicit basis $\{\Theta_i\}$ of the group algebra $K[G]$ over $K$,
which may depend upon the extension $L/K$ but should be independent of
the element $\rho$, with the additional property that the valuations
associated with this basis, $\{v_L(\Theta_i\rho)\}$, give a complete
set of residues modulo $[L:K]$.  These two ingredients, an integer
certificate and a basis, make up what we call a {\em Galois
scaffolding}.

\paragraph{\bf Prototype: Cyclic extensions of degree $p$} 
Let $L/K$ be a ramified, cyclic, degree $p$ extension of local fields with
$\mbox{Gal}(L/K)=\langle\sigma\rangle$.  Assume that the ramification
break number for $L/K$ is $b$ and $\gcd(p,b)=1$.  Note that this does
not restrict the extension when $K$ has characteristic $p$ and is only
a minor restriction when $K$ has characteristic $0$ \cite[III.~Prop
2.3]{fesenko}.  Let $\rho\in L$ be any element with $v_L(\rho)\equiv
b\bmod p$.  Then $v_L((\sigma-1)^i\rho)\equiv(i+1)b\bmod p$ for $0\leq
i\leq p-1$. In particular, $v_L((\sigma-1)^i\rho)$ yields a complete
set of residues modulo $p$, and so we have a Galois scaffolding: Pick
any integer $\equiv b\bmod p$ and the basis, $\{(\sigma-1)^i:0\leq i\leq
p-1\}$.

Galois scaffolding should be viewed as normal bases with the important
advantage that the valuation of any element expressed in terms of the
Galois scaffolding can be easily determined.  In the example above,
since $L/K$ is fully ramified, every element $\alpha\in L$ can be
expressed as $\alpha=\sum_{i=0}^{p-1}a_i(\sigma-1)^i\rho$ for certain
$a_i\in K$. Then $v_L(\alpha)=\min\{v_L(a_i)+ib+v_L(\rho):0\leq i\leq
p-1\}$.  We repeat ourselves for emphasis.  Normal bases and power
bases (polynomial bases) in a prime element are two common bases. The
first allows the Galois action to be easily followed. The second
allows for an easy determination of valuation.  These two properties are
usually at tension and so Galois scaffolding are remarkable for the
delicate balance that they achieve\footnote{It is easy to see that
Galois scaffolding are not universally available. Considering any
unramified extension, where there can be no integer certificate.}.

Galois scaffolding in ramified, cyclic, degree $p$ extensions have
made Galois module structure in these extensions tractable
\cite{ferton,borevich,aiba,desmit}, along with Galois module structure
in fully ramified, cyclic, degree $p^2$ extensions
\cite{elder:annals}.  In this paper, we will restrict our attention to
fully ramified elementary abelian extensions of local function fields
that are, in a particular sense, as simple as a ramified cyclic
extension of degree $p$, and give an explicit Galois scaffolding for
these extensions.  We are motivated by the fact that much about Galois
module structure in wildly ramified extensions remains poorly
understood despite the topic's venerable age.

\subsection{Notation}
Let $p$ be a prime integer and let $\mathbb{F}_p$ be the finite field
with $p$ elements.  Let $K=\mathbb{F}((t))$ be a local function field
with residue field $\mathbb{F}$, which is either $\mathbb{F}_q$, a
finite field with $q$ elements where $q$ is a power of $p$, or
$\overline{\mathbb{F}}_p$, the algebraic closure.  Let
$\wp:K\rightarrow K$ denote the $\mathbb{F}_p$-linear map
$\wp(x)=x^p-x$, and let $\phi$ denote the ring homomorphism
$\phi(x)=x^p$. Use subscripts to denote field of reference. So $\pi_K$
is a prime element of $K$, and $v_K$ is the valuation normalized so
that $v_K(\pi_K^t)=t$. Let $\euO_K=\{x\in K:v_K(x)\geq 0\}$ be the
valuation ring, and let $\euP_K=\pi_K\euO_K$ be its maximal ideal.
Let $L/K$ denote a fully ramified, Galois $p$-extension, with
$G=\mbox{Gal}(L/K)$. Define its ramification filtration by
$$G_i=\{\sigma\in G:v_L((\sigma-1)\pi_L)\geq i+1\}.$$

\subsection{One-dimensional elementary abelian extensions}
It is a basic observation in Artin-Schreier Theory that the elementary
abelian extensions of $K$ lie in one-to-one correspondence with the
finite subspaces of the $\bF_p$-vector space, $K/K^\wp $, where $K^\wp
$ denotes the image of $\wp $.

Assume for the moment that the residue field of $K$ is algebraically
closed, $\mathbb{F}=\overline{\mathbb{F}}_p$.  Define
$K_{(n)}=\phi^n(K)=\mathbb{F}((t^{p^n}))$ for $n\geq 1$.  Of course,
$K/K_{(n)}$ is an inseparable field extension, and so, in particular,
$K$ is a vector space over $K_{(n)}$.  Since the residue field of $K$
is algebraically closed, $K/K^\wp$ is also a vector space over
$K_{(n)}$.  We define {\em one-dimensional elementary abelian
extensions} to be those fully ramified elementary abelian extensions
of degree $p^i$ with $i\leq n+1$ that correspond to an $i$-dimensional
$\bF_p$-subspace of a one-dimensional $K_{(n)}$-subspace of
$K/K^\wp$. Of course, we are principally interested in maximal
extensions where $i=n+1$.

More generally, we can include the finite residue field case and define
one-dimensional elementary abelian extensions of degree $p^{n+1}$ to
be those that can be expressed as $L=K(x_0,\ldots ,x_n)$ with
$\wp(x_i)=x_i^p-x_i=\phi^n(\Omega_i)\cdot \beta$ for some $\beta\in
K$ with $v_K(\beta)=-b$, $b>0$ and $\gcd(b,p)=1$; and some $\Omega_i\in K$
that span an $n+1$-dimensional subspace over $\bF_p$, with $\Omega_0=1$ and 
$$v_K(\Omega_n)\leq \cdots \leq v_K(\Omega_1)\leq v_K(\Omega_0)=0.$$
Without any loss of generality, we can assume  moreover that whenever
$v_K(\Omega_i)=\cdots =v_K(\Omega_j)$ for $i<j$, the projections
of $\Omega_i, \ldots \Omega_j$ into
$\phi^n(\Omega_i)\beta\euO_K/\phi^n(\Omega_i)\beta\euP_K$ are linearly
independent over $\mathbb{F}_p$.  It should be clear from this
construction that the upper ramification numbers in one-dimensional
elementary abelian extensions of degree $p^{n+1}$ are congruent to
each other modulo $p^n$. Of course, the converse is not necessarily
true.

Simple examples of a one-dimensional elementary abelian extensions are
\begin{enumerate}
\item[(1)] extensions of the form $K(y)$ with $y^q-y=\beta$ (Lemma 5.2).
\end{enumerate}

It is probably not surprising that we are able to find Galois
scaffolding for a slightly broader class of extension, {\em near
one-dimensional elementary abelian extensions}, which arise when we
allow some error into the equations $\wp(x_i)=\phi^n(\Omega_i)\cdot
\beta$ defined above. In particular, we may replace those equations
with $\wp(x_i)\equiv \phi^n(\Omega_i)\cdot \beta+\epsilon_i$ for some
error terms $\epsilon_i\in K$ that satisfy a technical bound (6) and
use the same Galois scaffolding as for one-dimensional elementary
abelian extensions.

Simple examples of near one-dimensional elementary abelian extensions are
\begin{enumerate}
\item[(2)] fully ramified biquadratic extensions (Lemma 5.1), and
\item[(3)] fully and weakly ramified $p$-extensions (Lemma 5.3).
\end{enumerate}

Evidently, our Galois
scaffolding is not effected by small errors.  This last observation
can be rephrased in terms of twists by characters of Galois
representations, along the lines of \cite[\S2.2.3]{elder:necbreaks}.

\subsection{Galois scaffolding}
Assume the notation of the previous section and assume that $L/K$ is
near one-dimensional elementary abelian.

Relabel $\Omega_j^{(0)}=\Omega_j$, and perform the following
elementary row operations on the matrix
$[\phi^i(\Omega_j^{(0)})]_{0\leq i,j\leq n}$, which in passing we note
resembles the square root of a discriminant matrix.  The first column
is a column of $1$'s.  So start with the $i=n$ row and work down to
the $i=1$ row, subtracting the $i-1$st row from the $i$th row.  The
$i=0$ row and $i=0$ column of our matrix now agree with (1) below.  If
we ignore them, the result is
$[\phi^{i-1}(\wp(\Omega_j^{(0)}))]_{1\leq i,j\leq n}$. Divide each
entry in a row by the first entry of the row.  The result is
$[\phi^{i-1}(\wp(\Omega_j^{(0)})/\wp(\Omega_1^{(0)}))]_{1\leq i,j\leq
n}$.  Define $\Omega_j^{(1)}=\wp(\Omega_j^{(0)})/\wp(\Omega_1^{(0)})$
for $1\leq j\leq n$. Observe that $v_K(\Omega_n^{(1)})\leq \cdots \leq
v_K(\Omega_1^{(1)})=0$ and that the $\{\Omega_j^{(1)}\}_{1\leq j\leq
n}$ span an $n$ dimensional vector space over $\mathbb{F}_p$. Again we
have a matrix $[\phi^{i-1}(\wp(\Omega_j^{(1)}))]_{1\leq i,j\leq n}$
whose first column is a column of $1$'s.  Again, starting with the
$i=n$ row and working down to the $i=2$ row, we subtract the $i-1$st
row from the $i$th row. If we continue, following the same sequence of
steps as above, and repeat as often as necessary, we get
\begin{equation}
[\Omega]=\begin{bmatrix}
1 & \Omega_1^{(0)} & \Omega_2^{(0)} &\cdots &\Omega_n^{(0)} \\ 
0 & 1 &\Omega_2^{(1)} &\cdots &\Omega_n^{(1)} \\ 
& & \ddots & & &\\
0 & 0 & \cdots &1& \Omega_n^{(n-1)}\\
0 & 0 & \cdots &0& 1
\end{bmatrix}
\end{equation}
where $\Omega_j^{(0)}=\Omega_j$ and $\Omega_j^{(j)}=1$ for $0\leq
j\leq n$; and the $\Omega_j^{(i)}\in K$ for $1\leq i\leq n$ and $j> i$
are defined recursively by
$\Omega_j^{(i)}=\wp(\Omega_j^{(i-1)})/\wp(\Omega_{i}^{(i-1)})$.
Apply $\phi^{n-i-1}$ to row $i$ of $[\Omega]$, and get
$$[\Omega^\phi]=[\phi^{n-i-1}(\Omega_j^{(i)})]_{0\leq i,j\leq n}.$$

If we define the binomial coefficient $\binom{X}{i}$ by
$X\cdot(X-1)\cdots(X-i+1)/i!\in\mathbb{Q}[X]$, then we can define {\em
truncated exponentiation} to be the polynomial that results from the
truncation of the binomial series at the $p$th term:
\begin{equation}
(1+X)^{[Y]}:=\sum_{i=0}^{p-1}\binom{Y}{i}X^i\in\mathbb{Z}_{(p)}[X,Y]
\end{equation}
where $\mathbb{Z}_{(p)}$ denotes the integers localized at $p$. 

Choose $\sigma_i\in G=\mbox{Gal}(L/K)$ based upon our choice of
generators for $L/K$ by asking that
$$[(\sigma_i-1)x_j]=[\delta_{ij}]=I$$ ({\em i.e.}
$\sigma_ix_i=x_i+1$ and $\sigma_ix_j=x_j$ for $j\neq i$).  Define
$$[\Delta_{i,j}]=[\Omega^\phi]^{-1}.$$
Now
for $0\leq i\leq n$ define $\Theta_{(i)}\in
K[\sigma_n,\sigma_{n-1},\ldots ,\sigma_{n-i}]$ recursively by
\begin{equation}
\Theta_{(i)}=\sigma_{n-i}\Theta_{(0)}^{[-\Delta_{n-i,n}]}
\Theta_{(1)}^{[-\Delta_{n-i,n-1}]}\cdots
\Theta_{(i-1)}^{[-\Delta_{n-i,n-(i-1)}]}.
\end{equation}

Note that each
$\Theta_{(i)}$ is a 1-unit, {\em i.e.}  $\Theta_{(i)}\in
1+(\sigma-1:\sigma\in G)\subseteq K[G]$ where
$(\sigma-1:\sigma\in G)$ can be viewed as the augmentation ideal, the
Jacobson radical, or the nilradical. In particular, $\alpha^p=0$ for all
$\alpha\in (\sigma-1:\sigma\in G)$.  This means that
$(\Theta_{(i)}-1)^p=0$, and so
$\Theta_{(i)}^{[\Delta_{j,k}]}\Theta_{(i)}^{[-\Delta_{j,k}]}=1$.  As a
result, and since $\Delta_{n-r,n-r}=1$, the recursive definition for
$\Theta_{(i)}$ can be rewritten as
$$\sigma_{n-i}
=\Theta_{(0)}^{[\Delta_{n-i,n}]}\Theta_{(1)}^{[\Delta_{n-i,n-1}]}\cdots
\Theta_{(i-1)}^{[\Delta_{n-i,n-(i-1)}]}\Theta_{(i)}^{[\Delta_{n-i,n-i}]},$$
which
suggests 
the matrix equation:
$$\begin{bmatrix}
\Delta_{0,0} & \Delta_{0,1} &\cdots & \Delta_{0,n} \\ 
0 & \Delta_{1,1} &\cdots & \Delta_{1,n} \\ 
& & \ddots &  &\\ 
0 &\cdots &0 & \Delta_{n,n}
\end{bmatrix}\cdot \begin{bmatrix}
\Theta_{(n)}\\\Theta_{(n-1)}\\\vdots
\\\Theta_{(0)}\end{bmatrix}=\begin{bmatrix} \sigma_0\\\sigma_1\\\vdots
\\\sigma_n\end{bmatrix} ,$$ where addition is replaced by
multiplication and scalar multiplication is replaced truncated
exponentiation. Since truncated exponentiation does not distribute,
$(\Theta_{(i)}\Theta_{(j)})^{[\Delta]}\neq \Theta_{(i)}^{[\Delta]}
\Theta_{(j)}^{[\Delta]}$ (which is easy to check with $p=2$), we have
$[\Theta_{(n-j)}]\neq [\Omega^\phi]\cdot[\sigma_i]$,
despite the fact that $[\Omega^\phi]=[\Delta_{i,j}]^{-1}$. In other
words, this matrix equation is simply a convenient way to express a
recursive definition -- no more, no less.

We are prepared to state the main result of the paper, which is proven
in \S3, \S4.

\begin{theorem} 
Let $L/K$ be a near one-dimensional elementary abelian extension, as defined in \S1.2.  Let
$\Theta_{(i)}\in K[\mbox{\rm Gal}(L/K)]$ be defined as in {\rm (3)}.
For $1\leq i\leq n$, let $m_i=v_K(\Omega_{i-1})-v_K(\Omega_i)$, and
choose any $\alpha_j\in K$ with
$v_K(\alpha_j)=p^{n-j-1}\sum_{i=j+1}^{n}p^im_i$. Let $b_m$ be the
largest {\rm (}lower{\rm )} ramification break number of $L/K$.
Given any $\rho\in L$ with $v_L(\rho)\equiv b_m\bmod p^{n+1}$ and
$a_s\in\{0,\ldots , p-1\}$,
$$v_L\left (\prod_{s=0}^{n}\alpha_{n-s}^{a_s}(\Theta_{(s)}-1)^{a_s}\rho\right
)=v_L(\rho)+\sum_{s=0}^{n} a_sp^sb_m.
$$
\end{theorem}

As the integers $\sum_{s=0}^{n} a_sp^s$ run
through all possibilities from $0$ to $p^{n+1}-1$, the integers
$(\sum_{s=0}^{n} a_sp^s)b_m$ run through all residues modulo $p^{n+1}$.
Therefore
\begin{corollary} 
$L$ has a Galois scaffolding.
\end{corollary} 
\begin{corollary} 
Any element in $L$ of valuation $b_m$ generates a normal field basis.
\end{corollary} 
This last corollary provides evidence for the Conjecture in
\cite{elder:blms}.

\section{Cyclic extensions of degree $p$}

This paper is concerned with Galois, fully and thus wildly ramified
$p$-extensions that are, in a certain sense, as simple as cyclic
extensions of degree $p$.  And so, we should take a moment to consider
the prototype: Let $L/K$ be a cyclic, ramified extension of degree
$p$.  So $L=K(x)$ where $x$ satisfies $\wp(x)=x^p-x=\beta$ for some
$\beta\in K$ with $v_K(\beta)=-b$, $b>0$ and $\gcd(b,p)=1$.  Let
$\langle\sigma\rangle=\mbox{Gal}(L/K)$ with $\sigma x=x+1$. Since
$v_L((\sigma-1)x)=0$, it is easy to see that the integer $b$ is the
ramification break number for $L/K$.  Since $\wp(x)=\beta$ is really a
statement about the norm of $x$, namely $N_{L/K}(x)=\beta$, we have
$v_L(x)=-b$ as well.

We may rewrite $x^p-x=\beta$ as $x\cdot
\binom{x-1}{p-1}=-\beta$, where $\binom{x-1}{p-1}$ is a binomial
coefficient. Then 
$$\binom{x-1}{p-1}\in L$$ 
generates $L/K$, satisfies
$v_L(\binom{x-1}{p-1})=-(p-1)b\equiv b\bmod p$ and, we contend,
is a particularly insightful element 
to consider.  Recall the definition of truncated
exponentiation and notice the striking similarity between
$$\sigma^{[i]}\binom{x-1}{p-1}=\sigma^i\binom{x-1}{p-1}= \binom{x-1+i}{p-1} \mbox{ for } 0\leq i\leq p-1,$$
and the equation in
\begin{lemma} Let $L=K(x)$ with $x^p-x=\beta\in K$ be a cyclic extension 
with $\langle\sigma\rangle=\mbox{\rm Gal}(L/K)$. 
Given $A\in L$, 
$$\sigma^{[A]}\binom{x-1}{p-1}=\binom{x-1+A}{p-1}.$$
\end{lemma}
\begin{proof}
Recall Pascal's identity
$\binom{X}{i-1}+\binom{X}{i}=\binom{X+1}{i}\in \mathbb{Q}[X]$, which
can be rewritten as $\binom{X+1}{i}-\binom{X}{i}=\binom{X}{i-1}$. This
leads to the  nice observation, used in \cite{desmit}, that
$(\sigma-1)\binom{x-1}{i}=\binom{x-1}{i-1}$ for $0\leq
i\leq p-1$, and therefore
$$(\sigma-1)^i\binom{x-1}{p-1}=\binom{x-1}{p-1-i}\quad\mbox{for }0\leq
i\leq p-1.$$
 Under the substitution $X=\sigma-1$ and
$Y=A\in L$, we find
$$\sigma^{[A]}\binom{x-1}{p-1}=
\sum_{i=0}^{p-1}\binom{A}{i}(\sigma-1)^i\binom{x-1}{p-1}=
\sum_{i=0}^{p-1}\binom{A}{i}\binom{x-1}{p-1-i}\in L.$$ Vandermonde's
Convolution Identity $\sum_{i=0}^{p-1}\binom{X}{i}\binom{Y}{p-1-i}=
\binom{X+Y}{p-1}\in\mathbb{Z}_{(p)}[X,Y]$ results from considering the
coefficient of $Z^{p-1}$ in the identity
$(1+Z)^X(1+Z)^{Y}=(1+Z)^{X+Y}\in\mathbb{Q}[X,Y][[Z]]$. If we replace
$X=A$ and $Y=x$, we find
$\sum_{i=0}^{p-1}\binom{A}{i}\binom{x-1}{p-1-i}= \binom{x-1+A}{p-1}\in
L$.\end{proof} 

In \cite{elder:newbreaks} a refined ramification filtration was
introduced. It grew out of the possibility that the natural
$\mathbb{F}_p$-action on $\sigma$ could be extended to a residue field
``action;'' a possibility that is certainly suggested by this striking
similarity.

In this paper, we will develop a Galois scaffolding based on this similarity.
Specifically, we suppose that $L/K$ sits in a more general Galois
extension $M/N$, and we suppose furthermore that $L/N$ is normal and
that $\gamma\in\mbox{Gal}(M/N)$. So $\gamma^{-1} x=x+\delta$ for some
$\delta\in L$ and
$\sigma^{[\delta]}\binom{x-1}{p-1}=\binom{x-1+\delta}{p-1}=\gamma^{-1}\binom{x-1}{p-1}$. But
then
$$\gamma\sigma^{[\delta]}\binom{x-1}{p-1}=\binom{x-1}{p-1}.$$ If
$\delta\neq 0$ and $\gamma\not\in\langle\sigma\rangle$, then neither
$\sigma$ nor $\gamma$ individually fix the field generator
$\binom{x-1}{p-1}$. Yet together, using truncated exponentiation, they
do. As a result, if we suppose that $\delta\in N$ then the stabilizer
of $\binom{x-1}{p-1}$ in $N[\mbox{Gal}(M/N)]$ is larger than expected.

\section{Galois scaffolding}

This section is motivated by the observation of \S2 concerning the
stabilizer of $\binom{x-1}{p-1}$ and should be considered ``top-down.''
We begin with an generic abelian $p$-extension, which we ``organize''
using the ramification filtration. This ``organization'' defines a
matrix $[\Delta]$.  If the coefficients of $[\Delta]$ lie in our base
field $K$, the extension satisfies a strong assumption, which makes it
possible for us to construct a Galois scaffolding, but also makes the
extension elementary abelian.  At the end of the section, one question
remains: Are there any elementary abelian extensions that satisfy this
strong assumption?  In \S4 we construct extensions that do -- from the
``bottom-up.''

Let $K_n/K$ be a fully
ramified, abelian extension of degree $p^{n+1}$. The case $n=0$
was addressed in \S1. So assume $n\geq 1$.  Let $G=\mbox{Gal}(K_n/K)$
and let $G_i=\{\sigma\in G: v_n((\sigma-1)\pi_n)\geq i+1\}$ denote the
Hilbert ramification groups with break numbers $b_1<b_2<\cdots<b_m$
such that $G=G_{b_1}$, $G_{b_i}\supsetneq G_{b_i+1}=G_{b_{i+1}}$ and
$G_{b_m+1}=\langle e\rangle$. Because $K$ is characteristic $p$,
$\gcd(b_1,p)=1$, and by \cite[IV\S2 Prop 11]{serre:local}, $b_i\equiv
b_1\bmod p$.

Organize the extension by choosing a 
filtration of $n+1$ subgroups that include the Hilbert ramification
groups and satisfy $G_{(i)}/G_{(i+1)}\cong C_p$,
$$G=G_{(0)}\supsetneq G_{(1)}\supsetneq \cdots\supsetneq
G_{(n)}\supsetneq G_{(n+1)}=\langle e\rangle.$$ Indeed, since each
quotient of consecutive Hilbert ramification groups is elementary
abelian, this is easy to do. The result is a set $\{\sigma_0,\sigma_1,
\ldots ,\sigma_n\}$ that generates $G$ (though probably not a minimal
generating set), such that $G_{(i)}=\langle \sigma_i,\sigma_{i+1},
\ldots , \sigma_n\rangle$ and the projection of $\sigma_i$ generates
$G_{(i)}/G_{(i+1)}\cong C_p$.  For $i\geq 0$, let the fixed field
of $G_{(i)}$ be $K_{i-1}$, with $K_{-1}=K$ and define
$b_{(i)}=v_n((\sigma_i-1)\pi_n)-1$. This means that $b_{(0)}\leq
b_{(1)}\leq \ldots \leq b_{(n)}$ is a list of $n+1$ not necessarily
distinct integers and $\{b_{(0)},\ldots
,b_{(n)}\}=\{b_1,\ldots , b_m\}$.

Since $K_n/K$ is abelian, the Theorem of Hasse-Arf states that the
upper ramification numbers are integers \cite[IV\S3]{serre:local},
which is equivalent to
$b_i\equiv b_m\bmod [G:G_{b_{i+1}}]$ for $1\leq i\leq m$, and also 
to
\begin{equation}
b_{(i)}\equiv b_{(n)}\bmod p^{i+1}\mbox{ for }0\leq i\leq n.
\end{equation}

Since $\{b_{(0)},\ldots ,b_{(n)}\}$ is the set of ramification break
numbers for $K_n/K$, the ramification break numbers for $K_i/K$ are
$\{b_{(0)},\ldots ,b_{(i)}\}$ \cite[IV \S1 Prop 3 Cor]{serre:local}.
Altogether,
$\mbox{Gal}(K_i/K_{i-1})=G_{b_{(i)}}/G_{b_{(i+1)}}=\langle\bar{\sigma}_i\rangle\cong
C_p$, with $K_i/K_{i-1}$ having ramification break number $b_{(i)}$.
As a result, there are $X_i\in K_i$ such that $v_i(X_i)=-b_{(i)}$,
$\wp(X_i)=X_i^p-X_i=B_i\in K_{i-1}$ and $\sigma_iX_i=X_i+1$.  Define
$$\Delta_{i,j}=(\sigma_i-1)X_j.$$ So $\Delta_{i,j}=0$ when $i>j$,
and 
$\Delta_{i,i}=1$. Because 
$X_j\in K_j$ and
$\sigma_i\sigma_j=\sigma_j\sigma_i$, we have
$\Delta_{i,j}\in K_{j-1}$ when $i<j$. Furthermore,
$v_j(\Delta_{i,j})=v_j((\sigma_i-1)X_j)=b_{(i)}-b_{(j)}\leq 0$.
Collect these $\Delta_{i,j}$ into a matrix, whose $j$th column lies in
$K_{j-1}$,
$$[\Delta]=\begin{bmatrix}
\Delta_{0,0} & \Delta_{0,1} &\cdots & \Delta_{0,n} \\ 
0 & \Delta_{1,1} &\cdots & \Delta_{1,n} \\ 
& & \ddots &  &\\ 
0 &\cdots &0 & \Delta_{n,n}
\end{bmatrix}.
$$ 
Motivated by the final comment in \S2, and the fact that we want 
a basis for $K[G]$ over $K$, we
impose
\begin{assumption}
$\Delta_{i,j}\in K\mbox{ for all }0\leq i,j\leq n$.
\end{assumption}
\begin{lemma}
Under Assumption  1,
$K_n/K$ is elementary abelian.
\end{lemma}
\begin{proof}
Since $\Delta_{i,j}\in K$, we have
$\sigma_i^kX_j=X_j+k\Delta_{i,j}$ for $0\leq k\leq p$. This means that
$\sigma_i^pX_j=X_j$ for all $0\leq i,j\leq n$, and in particular,
$\sigma_i^pX_n=X_n$ for all $0\leq i\leq n$.  Since
$v_n(X_n)=-b_{(n)}$, we have $\gcd(v_n(X_n),p)=1$ and thus
$K_n=K(X_n)$.\end{proof}

We will proceed in three steps towards our Galois scaffolding.  First we choose
a nice element $\mathbb{X}\in K_n$ with $v_n(\mathbb{X})=b_{(n)}=b_m$.
Then we determine a basis for $K[G]$ over $K$ so that the valuations of
these basis elements applied to $\mathbb{X}$ yield a complete set of
residues mod $p^{n+1}$.  Finally we prove in Proposition 3.3 that this second step
continues to hold when $\mathbb{X}$ is replaced by any element of
valuation $b_m\bmod p^{n+1}$.

Define
$$\rho=\prod_{j=0}^n\binom{X_j}{p-1}\in K_n.$$ 
Because of (4), we may
choose $\alpha_j\in K$ such that $v_j(\alpha_j)=b_{(n)}-b_{(j)}$.
Therefore $v_j(\alpha_j^{-(p-1)}\binom{X_j}{p-1})=-(p-1)b_{(n)}$ for
$0\leq j\leq n$. Choose $\alpha\in K$ with $v_K(\alpha)=b_{(n)}$
Define $\mathcal{A}=\alpha\prod_{j=0}^n\alpha_j^{-(p-1)}\in K$. So
$v_n(\mathcal{A})\equiv 0\bmod p^{n+1}$ and
$$\mathbb{X}=
\mathcal{A}\rho
=
\alpha\prod_{j=0}^n\alpha_j^{-(p-1)}\binom{X_j}{p-1}$$
has valuation $v_n(\mathbb{X})=p^{n+1}b_{(n)}-(p-1)\sum_{j=0}^np^{n-j}b_{(n)}=b_{(n)}=b_m$.

Recall (3), namely the recursive definition for  $\Theta_{(i)}\in
K[G]$ for $0\leq i\leq n$.

\begin{lemma} For $0\leq i, j\leq n$,
$$
\Theta_{(i)}\binom{X_j}{p-1}=\begin{cases}
\binom{X_j}{p-1}&\mbox{if }j\neq n-i,\\
\binom{X_j+1}{p-1}&\mbox{if }j= n-i.
\end{cases}
$$
\end{lemma}
\begin{proof}
We proceed by induction.  For $i=0$,
$\Theta_{(i)}=\Theta_{(0)}=\sigma_n$ and since $\sigma_n$ fixes
$K_{n-1}$ while $\binom{X_j}{p-1}\in K_j$, the result is clear.  Now
assume the result for $0\leq i< k$ and consider
$\Theta_{(k)}\binom{X_j}{p-1}$. Because $\Theta_{(k)}$ is a product
(3), we need to examine the effect of each factor
$\Theta_{(i)}^{[-\Delta_{n-k,n-i}]}$ in that product, namely
$\Theta_{(i)}^{[-\Delta_{n-k,n-i}]}\binom{X_j}{p-1}$ for $0\leq i<k$.
By induction $(\Theta_{(i)}-1)^r
\binom{X_{n-i}}{p-1}=\binom{X_{n-i}}{p-1-r}$ for $0\leq r\leq p-1$,
and $(\Theta_{(i)}-1)^r \binom{X_j}{p-1}=0$ for $j\neq n-i$.
Therefore using Lemma 2.1, we have
$$\Theta_{(i)}^{[-\Delta_{n-k,n-i}]}\binom{X_j}{p-1}=\begin{cases}
\binom{X_j}{p-1}&\mbox{for }j\neq n-i,\\
\binom{X_j-\Delta_{n-k,j}}{p-1}&\mbox{for }j= n-i.
\end{cases}$$
If $j<n-k$, then every factor of $\Theta_{(k)}$ and thus
$\Theta_{(k)}$ acts trivially on $\binom{X_j}{p-1}$. If $j=n-k$ then
the only factor of $\Theta_{(k)}$ to act non-trivially is
$\sigma_{n-k}=\sigma_j$. As a result, $\Theta_{(k)}
\binom{X_j}{p-1}=\sigma_j\binom{X_j}{p-1}= \binom{X_j+1}{p-1}$. If
$j>n-k$, then exactly two factors of $\Theta_{(k)}$ to act
non-trivially, namely $\sigma_{n-k}$ and
$\Theta_{(n-j)}^{[-\Delta_{n-k,j}]}$. So
$$\Theta_{(k)}\binom{X_j}{p-1}=
\sigma_{n-k}\Theta_{(n-j)}^{[-\Delta_{n-k,j}]}\binom{X_j}{p-1}
=\sigma_{n-k}\binom{X_j-\Delta_{n-k,j}}{p-1}=\binom{X_j}{p-1}.$$\end{proof}

Now notice that for $0\leq r\leq p-1$, we have
\begin{multline*}(\Theta_{(i)}-1)^r\mathbb{X}=(\Theta_{(i)}-1)^r\mathcal{A}\prod_{j=0}^n\binom{X_j}{p-1}
=\mathcal{A}\prod_{j\neq i}\binom{X_j}{p-1}\cdot (\Theta_{(i)}-1)^r
\binom{X_{n-i}}{p-1}\\=\mathcal{A}\prod_{j\neq i}\binom{X_j}{p-1}\cdot 
\binom{X_{n-i}}{p-1-r}.
\end{multline*}
Therefore
$(\Theta_{(i)}-1)^r\mathbb{X}=\mathbb{X}\binom{X_{n-i}}{p-1-r}
\binom{X_{n-i}}{p-1}^{-1}$ and so
$v_n((\Theta_{(i)}-1)^r\mathbb{X})=b_n+rp^ib_{n-i}$. 
Moreover given
$c_i\in\{0,1\ldots ,p-1\}$, we have
$$\prod_{i=0}^n(\Theta_{(i)}-1)^{c_i}\mathbb{X}=\mathcal{A}\prod_{j=0}^n\binom{X_j}{p-1-c_{n-j}},$$
and
using the $\alpha_j\in K$ with
$v_j(\alpha_j)=b_{(n)}-b_{(j)}$,
\begin{equation}
v_n\left (\prod_{i=0}^n\alpha_{n-i}^{c_i}(\Theta_{(i)}-1)^{c_i}\mathbb{X}\right)=\left(1+\sum_{i=0}^nc_ip^i\right)b_{(n)}.\end{equation}
Therefore
$$\left \{\prod_{i=0}^n\alpha_{n-i}^{c_i}(\Theta_{(i)}-1)^{c_i}:0\leq c_i\leq p-1\right \}$$
is the desired basis.

\begin{proposition}
Under Assumption 1, we have a Galois scaffolding. Let $\mathbf{X}\in
K_n$ be any element with $v_n(\mathbf{X})\equiv b_{(n)}=b_m \bmod p^{n+1}$. Let
$\Theta_{(i)}\in K[G]$ be as defined in {\rm (3)}, and let
$\alpha_j\in K$ with $v_K(\alpha_j)=(b_{(n)}-b_{(j)})/p^{j+1}\in\mathbb{Z}$, then
$$v_n\left (\prod_{i=0}^n\alpha_{n-i}^{c_i}(\Theta_{(i)}-1)^{c_i}\mathbf{X}\right)=v_n(\mathbf{X})+\sum_{i=0}^nc_ip^ib_m.$$
\end{proposition}
\begin{proof}
Using (5), we can express $\mathbf{X}$ as a linear combination of
$\prod_{i=0}^n\alpha_{n-i}^{c_i}(\Theta_{(i)}-1)^{c_i}\mathbb{X}$ with
coefficients in $K$. It is enough therefore to show that 
when we apply $\prod_{i=0}^n\alpha_{n-i}^{d_i}(\Theta_{(i)}-1)^{d_i}$
with $0\leq
d_i\leq
p-1$ to any term in this linear combination, we increase valuation by at least
$\sum_{i=0}^nd_ip^ib_m$, namely that
$$v_n\left(\prod_{i=0}^n\alpha_{n-i}^{c_i+d_i}(\Theta_{(i)}-1)^{c_i+d_i}\mathbb{X}\right)\geq
v_n\left(\prod_{i=0}^n\alpha_{n-i}^{c_i}(\Theta_{(i)}-1)^{c_i}\mathbb{X}\right)+\sum_{i=0}^nd_ip^ib_m.$$
If any sum $c_i+d_i\geq p$ then $(\Theta_{(i)}-1)^{c_i+d_i}=0$ and the
valuation of the left-hand-side is infinite. So we are left with the
case where all sums $c_i+d_i<p$. But in this case, we can use (5) to
determine that we have equality.
\end{proof}

\section{Near One-dimensional Elementary Abelian Extensions}

In contrast with \S3, this section is ``bottom-up''. Motivated by the
idea of maximal refined ramification in \cite{elder:necbreaks}, we
follow \S1.2 and define the class of near one-dimensional elementary
abelian extensions, by describing how the generators of each extension
are related.  We organize these generators by size (by valuation) as
in \S1.2, and then define the matrix $[\Omega^\phi]$ over $K$ as in
\S1.3. Our organization of the
generators, ``organizes'' the matrix $[\Omega^\phi]$. The main result
of the section is that this also ``organizes'' the extension in
essentially the same fashion as in \S3. In particular,
$[\Omega^\phi]\cdot[\Delta]=I$, which means that near one-dimensional
elementary abelian extension satisfy Assumption 1 and thus possess
Galois scaffolding.

Recall the notation of \S1.2: Let $L=K(x_0, \ldots , x_n)$ with
$\wp(x_i)=\phi^n(\Omega_i)\cdot \beta +\epsilon_i$ for some $\beta\in
K$ with $v_K(\beta)=-b$, $b>0$ and $\gcd(b,p)=1$; some $\Omega_i\in K$
that span an $n+1$-dimensional subspace over $\bF_p$; and some ``error
terms'' $\epsilon_i\in K$, whose size will be controlled by (6)
below. Initially, we merely assume
$v_K(\epsilon_i)>v_K(\phi^n(\Omega_i)\beta)$, so the ramification
break number of $K(x_i)/K$ is $-v_K(\phi^n(\Omega_i)\beta)$.

Furthermore recall $\Omega_0=1$ and that the other $\Omega_i$ are
``organized'' (relabelled) so that $v_K(\Omega_n)\leq \cdots \leq
v_K(\Omega_1)\leq v_K(\Omega_0)=0$, and if $v_K(\Omega_i)=\cdots
=v_K(\Omega_j)$ for $i<j$, the projections of $\Omega_i, \ldots
\Omega_j$ into
$\phi^n(\Omega_i)\beta\euO_K/\phi^n(\Omega_i)\beta\euP_K$ are linearly
independent over $\mathbb{F}_p$.  This means that $K(x_i,\ldots ,x_j)$
has one break in its ramification filtration at
$-v_K(\phi^n(\Omega_i)\beta)$.  

For $1\leq i\leq n$, define
$m_i=v_K(\Omega_{i-1})-v_K(\Omega_i)\geq 0$.
We control the size of the error terms with: For
$1\leq i\leq n$,
\begin{multline}
v_K(\epsilon_i)>-\frac{b}{p^n}-\sum_{j=1}^ip^jm_j+
\sum_{j=i+1}^n(p^n-p^j)m_j\\
=v_K(\phi^n(\Omega_i)\beta)+\frac{(p^n-1)b}{p^n}-(p-1)\sum_{j=1}^{n-1}p^jv_K(\Omega_j),
\end{multline} 
which since $v_K(\Omega_j)\leq 0$ is clearly stronger than our initial
assumption, $v_K(\epsilon_i)>v_K(\phi^n(\Omega_i)\beta)$. Notice
further that if, for a particular $i$, the right-hand-side of (6) is
zero, then (6) is equivalent to ``no error'' ({\em i.e.}
$\epsilon_i=0$), since the inequality $v_K(\epsilon_i)>0$ implies
$\epsilon_i\in K^\wp$.

Choose $\sigma_i\in G=\mbox{Gal}(L/K)$ based upon our generators so
that $[(\sigma_i-1)x_j]=[\delta_{ij}]=I$.  Define $H_{(i)}=\langle
\sigma_i, \ldots , \sigma_n\rangle$, and let $K_{i-1}=K(x_0,\ldots
,x_{i-1})$ be the fixed field of $H_{(i)}$.  So $K_{-1}=K$ and
$K_n=L$. As noted earlier, 
$$u_{(i)}=b+p^n\sum_{j=1}^{i}m_j$$ is the ramification number of
$K(x_i)/K$, and is therefore an upper ramification number of $L/K$. By
considering our assumptions on the $\Omega_i$, one sees that the set
of upper ramification numbers is $\{u_{(0)}, \ldots , u_{(n)}\}$. We
may pass to the lower ramification numbers using the Herbrand function
$\psi(x)$ \cite[IV \S3]{serre:local}. Again considering our assumptions on
the $\Omega_i$, one sees that $\{b_{(0)}, \ldots ,b_{(n)}\}$ is the set of lower ramification numbers where
$$b_{(i)}=b+p^n\sum_{j=1}^{i}p^jm_j.$$ Moreover, $b_{(i)}$ is the
ramification number of $K_i/K_{i-1}$, and it is clear that the groups
$H_{(i)}$ are the groups $G_{(i)}$ defined in \S3. We can express the restriction on the error terms in (6) in terms of ramification numbers:
$v_K(\epsilon_i)>-b_{(n)}/p^n+u_{(n)}-u_{(i)}$.

Our next step is to construct the $X_i\in K_i$ of \S3. 
Recall
the $\Omega_j^{(i)}$ defined in \S1.3.
Define $X_j^{(0)}=x_j$.  And
for $j\geq i$, recursively define
\begin{equation}
X_j^{(i)}=X_j^{(i-1)}- \phi^{n-i}(\Omega_j^{(i-1)})
X_{i-1}^{(i-1)}.
\end{equation}
If we use this definition to replace $X_j^{(i-1)}$ in (7) with
$X_j^{(i-2)}- \phi^{n-i+1}(\Omega_j^{(i-2)}) X_{i-2}^{(i-2)}$, we find
that $X_j^{(i)}=X_j^{(i-2)}- \phi^{n-i+1}(\Omega_j^{(i-2)})
X_{i-2}^{(i-2)}- \phi^{n-i}(\Omega_j^{(i-1)}) X_{i-1}^{(i-1)}$.
If we continue in this way, we eventually find
$X_j^{(i)}=X_j^{(0)}-\sum_{k=0}^{i-1} \phi^{n-k-1}(\Omega_j^{(k)})
X_{k}^{(k)}$.  Consider the case $i=j$. Since $x_j=X_j^{(0)}$ and
$\Omega_j^{(j)}=1$, this can be rewritten as
$x_j=\sum_{k=0}^j\phi^{n-k-1}\Omega_j^{(k)}X_k^{(k)}$.  
Recall that
$[\Omega^\phi]=[\phi^{n-i-1}(\Omega_j^{(i)})]_{0\leq i,j\leq n}$.
Therefore
$$[X_0^{(0)},X_1^{(1)},\ldots ,X_n^{(n)}]
\cdot[\Omega^\phi]=[x_0,x_1,x_2,\cdots , x_n].$$ 
Since
$I=[(\sigma_i-1)x_j]$, 
we find that
$[(\sigma_i-1)X_j^{(j)}]\cdot[\Omega^\phi]=I$. Therefore
\begin{equation}
[(\sigma_i-1)X_j^{(j)}]=[\Omega^\phi]^{-1}.
\end{equation}

Clearly $K_j=K(x_0,\ldots ,x_j)=K(X_0^{(0)},\ldots ,X_j^{(j)})$.  If
we could determine that $v_j(X_j^{(j)})=-b_{(j)}$, then we could choose
the $X_j=X_j^{(j)}$ and find that $[\Delta_{i,j}]=[\Omega^\phi]^{-1}$.
As a result, our extension would satisfy Assumption 1.  The remainder
of this section is therefore concerned with the valuation
$v_j(X_j^{(j)})$.  Since the $\Omega_j^{(i)}$ are an important
ingredient in the definition of the $X_j^{(j)}$, given in (7), we need

\begin{lemma} For $0\leq i<j\leq n$
$$v_K(\Omega_j^{(i)})=-p^i\sum_{k=i+1}^jm_k$$
\end{lemma}
\begin{proof} We induct on $i$.
Since $m_k= v_K(\Omega_{k-1}^{(0)})-v_K(\Omega_k^{(0)})$ for $1\leq
k\leq n$, the result holds for $i=0$. 
For $i>1$, we assume the result. So in particular,
$v_K(\Omega_n^{(i-1)})\leq \cdots \leq v_K(\Omega_{i+1}^{(i-1)})\leq 
v_K(\Omega_i^{(i-1)})=0$. Then
$v_K(\wp(\Omega_j^{(i-1)}))= pv_K(\Omega_j^{(i-1)})$
and thus using the definition for $\Omega_j^{(i)}$ in \S1.3,
we find that
$v_K(\Omega_j^{(i)})= pv_K(\Omega_j^{(i-1)})
-pv_K(\Omega_i^{(i-1)}) $ and result follows.
\end{proof}

To assist in our analysis of $v_j(X_j^{(j)})$, define
$B_0=\beta$,
$E_j^{(0)}=\epsilon_j$ for $j>0$.
Then for $i>0$ recursively define
\begin{equation}
B_i=-\phi^{n-i}(\wp(\Omega_i^{(i-1)}))X_{i-1}^{(i-1)}+E_i^{(i-1)}
\end{equation}
and $E_j^{(i)}=E_j^{(i-1)}- \phi^{n-i}(\Omega_j^{(i)}) E_i^{(i-1)}$
for $j>i$.  And $E_i^{(i)}=0$. The significance of these $B_i$ and
$E_j^{(i)}$ results from
\begin{lemma} For $j\geq i$
$$\wp(X_j^{(i)})=\phi^{n-i}(\Omega_j^{(i)})B_i+E_j^{(i)}$$
\end{lemma}
\begin{proof}
The statement is clear for $i=0$. Assume that it holds for $i-1$.
Therefore 
$\wp(X_j^{(i-1)})=\phi^{n-i+1}(\Omega_j^{(i-1)})B_{i-1}+E_j^{(i-1)}$
and in particular, $\wp(X_{i-1}^{(i-1)})=B_{i-1}$.
Consider $\wp(X_j^{(i)})$.
It is easy to see that
$\wp(aX)=\phi(a)\wp(X)+\wp(a)X$. Therefore using (7) we find that
\begin{multline*}
\wp(X_j^{(i)})=\wp(X_j^{(i-1)})-
\phi^{n-i+1}(\Omega_j^{(i-1)})\wp(X_{i-1}^{(i-1)})-
\phi^{n-i}(\wp(\Omega_j^{(i-1)}))X_{i-1}^{(i-1)}\\
=\phi^{n-i+1}(\Omega_j^{(i-1)})B_{i-1}+ E_j^{(i-1)}-
\phi^{n-i+1}(\Omega_j^{(i-1)})B_{i-1}-
\phi^{n-i}(\wp(\Omega_j^{(i-1)}))X_{i-1}^{(i-1)}\\ = E_j^{(i-1)} -
\phi^{n-i}(\wp(\Omega_j^{(i-1)}))X_{i-1}^{(i-1)},
\end{multline*}
which, using (9), can be seen to agree with the statement for $i$.
\end{proof}

\begin{lemma} Assume the bounds given in {\rm (6)}. Then for
$1\leq i\leq n$, we have $$v_K(E_i^{(i-1)})>-b_{(i)}/p^i.$$
\end{lemma}
\begin{proof}
Use Lemma 4.1 to determine that (6) is equivalent to
$$v_K(\phi^{n-i}(\Omega_{n}^{(i)})\epsilon_i)>-b_{(n)}/p^n.$$
We are interested in $v_K(E_i^{(i-1)})$. So
recall that
$E_j^{(i)}=E_j^{(i-1)}- \phi^{n-i}(\Omega_j^{(i)}) E_i^{(i-1)}$
for $j>i$, which means that $E_j^{(i)}=E_j^{(0)}- \sum_{k=1}^{i}
\phi^{n-k}(\Omega_j^{(k)}) E_k^{(k-1)}$, and
in particular,
\begin{equation}
E_i^{(i-1)}=\epsilon_i- \sum_{k=1}^{i-1}
\phi^{n-k}(\Omega_{i}^{(k)}) E_k^{(k-1)}.
\end{equation}
In order that
$v_K(E_i^{(i-1)})>-b_{(i)}/p^i$ for $1\leq i\leq n$, it is sufficient to
prove 
\begin{eqnarray}
v_K(\epsilon_i)&>&-b_{(i)}/p^i\mbox{ for } 1\leq i\leq n,\mbox{ and}\\
v_K(\phi^{n-k}(\Omega_{i}^{(k)})
E_k^{(k-1)})&>&-b_{(i)}/p^i\mbox{ for } 1\leq k\leq i-1\leq n-1.
\end{eqnarray}

Let $A_i=-b_{(i)}/p^i+v_K(\phi^{n-i}(\Omega_n^{(i)}))$. Using Lemma
4.1, we find that
$-b_{(i)}/p^i+v_K(\phi^{n-i}(\Omega_n^{(i)}))=-b_{(i-1)}/p^i+v_K(\phi^{n-i+1}(\Omega_n^{(i-1)}))$. As
a result, $A_i>A_{i-1}$, since $-b_{(i-1)}/p^i>-b_{(i-1)}/p^{i-1}$. We
are given by (6) that
$v_K(\phi^{n-i}(\Omega_{n}^{(i)})\epsilon_i)>-b_{(n)}/p^n=A_n$.  So
$v_K(\phi^{n-i}(\Omega_{n}^{(i)})\epsilon_i)>A_j$ for all $j$, including
$j=i$. Therefore (11) follows from (6).

Focus on (12), which is equivalent to  
$v_K(E_k^{(k-1)})>B_i^k$ where
$B_i^k=
-b_{(i)}/p^i-v_K(\phi^{n-k}(\Omega_{i}^{(k)}))$. Since
$-b_{(i)}/p^i-v_K(\phi^{n-k}(\Omega_{i}^{(k)}))=
-b_{(i-1)}/p^i-v_K(\phi^{n-k}(\Omega_{i-1}^{(k)}))$, we have
$B_i^k>B_{i-1}^k$. And thus (12) is equivalent to
\begin{equation}
v_K(\phi^{n-k}(\Omega_{n}^{(k)})
E_k^{(k-1)})>-b_{(n)}/p^n\mbox{ for }1\leq k\leq n-1.
\end{equation} 

Switch the roles of $i$ and $k$ in (10) and then apply
$\phi^{n-k}(\Omega_{n}^{(k)})$ to both sides:
$$\phi^{n-k}(\Omega_{n}^{(k)})E_k^{(k-1)}=
\phi^{n-k}(\Omega_{n}^{(k)})\epsilon_k- \sum_{i=1}^{k-1}
\phi^{n-k}(\Omega_{n}^{(k)})\phi^{n-i}(\Omega_{k}^{(i)})
E_i^{(i-1)}.$$

By Lemma 4.1,
$v_K(\phi^{n-k}(\Omega_{n}^{(k)})\phi^{n-i}(\Omega_{k}^{(i)}))=
v_K(\phi^{n-i}(\Omega_{n}^{(i)}))$. Therefore (13) follows from (6) by
induction on $k$.
\end{proof}

\begin{lemma} Assume the bounds in {\rm (6)}. Then for $0\leq j\leq n$, $v_j(X_j^{(j)})=-b_{(j)}$.
\end{lemma}
\begin{proof}
It is clear that $v_0(X_{0}^{(0)})=-b_{(0)}$. So for $i>0$, assume
that
$v_{i-1}(X_{i-1}^{(i-1)})=-b_{(i-1)}=-b-p^n\sum_{j=1}^{i-1}p^jm_j$.
Using Lemma 4.1, we see that
$v_K(\wp(\Omega_i^{(i-1)}))=-p^im_i$. So
$v_K(\phi^{n-i}(\wp(\Omega_i^{(i-1)})))=-p^nm_i$ and therefore
$v_{i-1}(\phi^{n-i}(\wp(\Omega_i^{(i-1)})))=-p^n\cdot p^im_i$. So
$v_{i-1}(\phi^{n-i}(\wp(\Omega_i^{(i-1)})X_{i-1}^{(i-1)})=-b_{(i)}$.
By Lemma 4.3, $v_{i-1}(E_i^{(i-1)})>-b_{(i)}$. Therefore
$v_{i-1}(B_i)=-b_{(i)}$. 
Lemma 4.2 implies that in particular the norm
$N_{K_{i}/K_{i-1}}(X_i^{(i)})=\wp(X_i^{(i)})=B_i$, which means that
$v_{i}(X_i^{(i)})=-b_{(i)}$.
\end{proof}

As a result, we can put all this together and find
\begin{proposition}
Near one-dimensional elementary abelian extensions satisfy Assumption 1.
\end{proposition}

\section{Examples of near one-dimensional elementary abelian extensions}

\begin{lemma}
Fully ramified biquadratic extensions are near one-dimensional
elementary abelian extensions.
\end{lemma}
\begin{proof}
Biquadratic extensions are special in that there is only one
nontrivial residue modulo $2$. Let $L/K$ be a fully ramified
biquadratic extension. We may assume that $L=K(x_0,x_1)$ with
$x_0^2-x_0=\beta$, $x_1^2-x_1=\beta_1$, $v_K(\beta_1)\leq
v_K(\beta)<0$ and both of $v_K(\beta_1)$ and $v_K(\beta)$ odd.
Because the difference of two odd numbers is even, there is a
$\mu_0\in K$ such that $\mu_0^2\beta_1\equiv \beta\bmod
\beta\euP_K$. Let $\beta=\mu_0^2\beta_1+\tau_0$ for some
$v_K(\tau_0)>v_K(\beta)$.  Since we can replace $\beta$ by any element
in its coset $\beta+K^{\wp}$, we may assume $v_K(\tau_0)=0$, or
$v_K(\tau_0)<0$ with $v_K(\tau_0)$ odd.  If $v_K(\tau_0)$ odd, then
there is a $\mu_1\in K$ such that $\mu_1^2\beta_1\equiv \tau_0\bmod
\tau_0\euP_K$, and thus $\beta=(\mu_0+\mu_1)^2\beta_1+\tau_1$ for
$v_K(\tau_1)>v_K(\tau_0)$.  Continue in this way until
$\beta=\mu^2\beta_1+\tau$ for some $\mu\in K$ and either $\tau=0$ or
$v_K(\tau)=0$.

If $\tau=0$, then $\beta_1=\mu^{-2}\beta$ and the extension is
one-dimensional. If $v_K(\tau)=0$, then
$\beta_1=\Omega_1^2\beta+\epsilon_1$ where $\epsilon_1=-\tau\mu^{-2}$
and $\Omega_1=\mu^{-1}$. Continuing to translate into the notation of
\S4, we note that $b=-v_K(\beta)$ and
$m_1=-v_K(\Omega_1)=v_K(\mu)$. So $v_K(\epsilon_1)=-2m_1>-b/2-2m_1$,
which is the inequality given by (6). So the extension is near
one-dimensional.
\end{proof}
\begin{lemma}
Let $K=\mathbb{F}((t))$ with $\mathbb{F}_q\subseteq \mathbb{F}$, and
let $\beta\in K$ with $v_K(\beta)<0$ and $\gcd(v_K(\beta),p)=1$.  Then
$L=K(y)$ with $y^q-y=\beta$ is a one-dimensional
elementary abelian extension of $K$.
\end{lemma}
\begin{proof}
Let $q=p^f$ and let $\{1=\omega_0,\omega_1, \cdots ,\omega_{f-1}\}$ be
a basis for $\mathbb{F}_q$ over $\mathbb{F}_p$. Then
$x_i=\sum_{r=0}^{f-1}\phi^r(\omega_iy)$ where $y^q-y=\beta$ satisfies
$x_i^p-x_i=\omega_i\beta$. Of course $\phi$ is an automorphism of
$\mathbb{F}_q$. So we let may set $\Omega_i=\phi^{-f+1}(\omega_i)$.
\end{proof}

The following class of fully and weakly ramified $p$-extensions ({\em i.e.}
with $G=G_1$ and $G_2=\{e\}$) is notable for being wildly ramified while
possessing a normal integral basis (for the maximal ideal)
\cite{ullom}.

\begin{lemma}
Let $L/K$ be a noncyclic, fully and weakly ramified $p$-extension,
then $L/K$ is a near one-dimensional elementary abelian extension.
\end{lemma}
\begin{proof}
The extension is elementary abelian \cite[IV \S2]{serre:local}, with
one break in its ramification filtration at $b=1$. As a result there
is only one upper ramification break number, also at $u=1$. Thus
$L=K(x_0,x_1,\ldots x_n)$ with $v_K(\wp(x_i))=-1$. Let
$\beta=\wp(x_0)$. Then there are units $\omega_i\in\mathbb{F}$ such
that $\wp(x_i)=\omega_i\beta\bmod \euO_K$. Since $\phi$ is an
automorphism of $\mathbb{F}$, we may let
$\Omega_i=\phi^{-n}(\omega_i)$ and find $\epsilon_i\in \euO_K$ such
that $\wp(x_i)=\phi^n(\Omega_i)\beta+\epsilon_i$ with either
$\epsilon_i=0$ or $v_K(\epsilon_i)=0$. Using the notation of \S4, we
find that $m_i=0$ and in all cases $v_K(\epsilon_i)\geq
0>-1/p^n=-b/p^n$, which is (6).
\end{proof}

\bibliography{bib} \end{document}